\documentclass[conference]{IEEEtran}
\IEEEoverridecommandlockouts
\usepackage{amsmath}
\usepackage{amssymb}
\usepackage{amsfonts}
\usepackage{mathtools}
\usepackage{mathrsfs}
\usepackage{dblfloatfix}
\usepackage{graphicx}
\usepackage[x11names]{xcolor}
\usepackage{tikz}
\usepackage{pgfplots}
\pgfplotsset{compat=1.18}
\usepackage{svg}
\usepackage{accents}
\usepackage{subcaption}
\usepackage{siunitx}
\usepackage{multirow}
\usepackage{algorithm}
\usepackage{algpseudocode}

\usepackage{circuitikz}
\definecolor{color_WF}{rgb}{0.2,0.1333,0.5333}
\definecolor{color_GFPP}{rgb}{0.800,0.1600,0.0337}
\definecolor{color_blue_1}{rgb}{0.267004,0.004874,0.329415}
\definecolor{color_blue_2}{rgb}{0.229739,0.322361,0.545706}
\definecolor{color_blue_3}{rgb}{0.128729,0.563265,0.551229}
\definecolor{color_blue_4}{rgb}{0.360741,0.785964,0.387814}
\definecolor{color_blue_5}{rgb}{0.993248,0.906157,0.143936}
\ctikzset{bipoles/length=0.85cm}

\definecolor{IEE_light_blue}{HTML}{1E90FF}
\definecolor{IEE_blue}{HTML}{206173}
\definecolor{IEE_red}{HTML}{F70146}
\definecolor{IEE_green}{HTML}{78BE73}
\definecolor{IEE_orange}{HTML}{D58E00}

\newcommand*{\tran}{^{\mkern-1.5mu\mathsf{T}}}

\definecolor{red1}{HTML}{fdccda}
\definecolor{red2}{HTML}{f70146}

\usepackage{hyperref}
\usepackage{csquotes}
\usepackage[utf8]{inputenc}
\usepackage[
backend = biber,
style=ieee, 
sorting=none, 
natbib=true, 
hyperref=true, 
date=year,
isbn=false,
maxbibnames = 3, 
dashed=false, 
]{biblatex}
\AtEveryBibitem{
    \clearfield{urlyear}
    \iffieldundef{pages}{}{\clearfield{doi}}
    \ifentrytype{misc}{\iflistundef{publisher}{}{\clearfield{url}}}{\clearfield{url}}
}

\usepackage[english]{babel}
\usepackage{stackengine}

\newcommand\munderbar[1]{%
  \underaccent{\bar}{#1}}

\definecolor{dkred}{rgb}{0.8,0,0}
\definecolor{blue}{rgb}{0,0,1}

\newcommand{\NCP}{\mathbf{NCP}}
\newcommand{\phiUpVec}{\bar{\boldsymbol{\varphi}}}
\newcommand{\phiLoVec}{\munderbar{\boldsymbol{\varphi}}}

\newcommand{\pfVecOp}{\boldsymbol{\hat{f}}^{*}}

\def\BibTeX{{\rm B\kern-.05em{\sc i\kern-.025em b}\kern-.08em
    T\kern-.1667em\lower.7ex\hbox{E}\kern-.125emX}}
    
\begin{document}

\title{Congestion-Sensitive Grid Aggregation for DC Optimal Power Flow \\
\thanks{Funded by the European Union (ERC, NetZero-Opt, 101116212). Views and opinions expressed are however those of the author(s) only and do not necessarily reflect those of the European Union or the European Research Council. Neither the European Union nor the granting authority can be held responsible for them.}
}

\author{\IEEEauthorblockN{Benjamin Stöckl, Yannick Werner, and Sonja Wogrin}
\IEEEauthorblockA{Institute of Electricity Economics and Energy Innovation, Graz University of Technology, Graz, Austria \\
\{benjamin.stoeckl,yannick.werner,wogrin\}@tugraz.at}}

\maketitle

\begin{abstract}

    The vast spatial dimension of modern interconnected electricity grids challenges the tractability of the DC optimal power flow problem. Grid aggregation methods try to overcome this challenge by reducing the number of network elements. Many existing methods use Locational Marginal Prices as a distance metric to cluster nodes. In this paper, we show that prevalent methods adopting this distance metric fail to adequately capture the impact of individual lines when there is more than one line congested. This leads to suboptimal outcomes for the optimization of the aggregated model. To overcome those issues, we propose two methods based on the novel Network Congestion Price metric, which preserves the impact of nodal power injections on individual line congestions. The proposed methods are compared to several existing aggregation methods based on Locational Marginal Prices. We demonstrate all methods on adapted versions of the IEEE RTS 24- and 300-Bus systems. We show that the proposed methods outperform existing approaches both in terms of objective function value error and maximum line limit violation, while exhibiting faster node clustering. We conclude that aggregation methods based on the novel Network Congestion Price metric are better at preserving the essential physical characteristics of the network topology in the grid aggregation process than methods based on Locational Marginal Prices. 

\end{abstract}


\section{Introduction}

    Achieving a carbon net-zero energy system requires increasing electricity generation from renewable energy sources~(RES). Large and interconnected electricity grids can facilitate the efficient integration of RES by balancing fluctuations~\cite{allard_european_2020}. System operators widely use Power System Optimization Models (PSOMs) to ensure a cost-efficient operation of large electricity grids with many sources of intermittent RES~\cite{hoffmann_review_2024}. They often employ linear DC power flow approximations to simplify the nonlinear and nonconvex AC power flow equations, which govern the flow of electricity in physical power grids~\cite{molzahn_survey_2019}. The vast spatial complexity of modern interconnected electricity grids, however, challenges the tractability of PSOMs despite employing a DC optimal power flow~(DC-OPF). By reducing the spatial dimension of PSOMs, grid aggregation methods frequently address this challenge~\cite{hoffmann_review_2024}. These methods try to conserve the essential physical properties of a given network topology with fewer network elements~\cite{fortenbacher_transmission_2018}.

    The aggregation process can be classified into two groups: spatial aggregation (e.g., nodes and lines) and technological aggregation (e.g., generators and demands). Spatial aggregation can be further divided into identifying groups of nodes, which we refer to as \textit{grid partitioning}, and determining parameters of the lines in the aggregated grid, which we refer to as \textit{grid aggregation}. In this paper, we neglect technological aggregation and focus on grid partitioning only, enabling an unbiased and rigorous comparison of different methods. We refer the interested reader to~\cite{hoffmann_review_2024} for further information on technological aggregation and provide a brief summary of existing grid aggregation methods in  Appendix~\ref{s:appendix_grid_aggregation}.
    
    In recent years, several methods combining various distance metrics with clustering algorithms to aggregate electrical grids have been proposed in the literature~\cite{hoffmann_review_2024, chicco_overview_2019}. In the context of DC-OPF, prevalent distance metrics based solely on the network topology include the geographical~\cite{phillips_spatial_2023} or electrical distance~\cite{cotilla-sanchez_multi-attribute_2013, biener_grid_2020}, or capacities of renewable energy sources~\cite{akdemir_open-source_2024}. These distance metrics do generally not consider specific operational conditions or market outcomes. While this makes them readily applicable, it may lead to suboptimal model outcomes for the aggregated grid. To harness this untapped potential, several studies~\cite{singh_reduced_2005, imran_effectiveness_2008, cao_incorporating_2018} have proposed leveraging market outcomes, mainly Locational Marginal Prices (LMPs), during the aggregation. The idea is to identify nodes with similar LMPs and group them together. This has also been the guiding principle in the determination of bidding zones in zonal markets, such as in the recent bidding zone review in Europe~\cite{acer_acer_decision_17-2023_baltic_bz_configurations_nodate}.

    Using the LMPs as a distance metric without additional information on the grid topology, however, has a major drawback when multiple line congestions exist, as shown in~\cite{jakubek_are_2015, colella_model-based_2021}. As we revisit in this paper, the impact of individual line congestions gets aggregated in the LMPs. Consequently, utilizing the LMPs as a distance metric may produce grid aggregations that fail to preserve individual line congestions, leading to infeasible operational schedules in the full model. To overcome this problem, we propose the novel Network Congestion Price (NCP) distance metric, which we combine with two clustering algorithms. Inspired by a time series aggregation method based on active constraints~\cite{wogrin_time_2023}, the NCP leverages the results of a DC-OPF problem while preserving the impact of individual line congestions. We demonstrate that this metric can generate partitions for grid aggregation sensitive toward multiple congestions, while the aggregated model obtains minimal error in the objective function value.
    
    The remainder of this paper is structured as follows. Section~\ref{s:model_formulation} briefly recaps the DC-OPF problem. Afterward, Section~\ref{s:appendix_grid_aggregation_methodology} recalls the derivation of the LMPs and introduces the novel NCP metric. Furthermore, it presents five distinct methods for grid partitioning for DC-OPF and a process to evaluate their quality. Section~\ref{s:case_study_results} demonstrates and compares all grid aggregation methods on two stylized case studies, followed by a conclusion in Section~\ref{s:conclusion}.

\section{DC optimal power flow problem}\label{s:model_formulation}

    Before recapping the lossless DC-OPF problem, we introduce the sets, variables, and parameters.
    Let $n \in \mathcal{N}$ be the set of buses with $|\mathcal{N}| = N$, $l \in \mathcal{L}$ be the set of lines with $|\mathcal{L}| = L$ and $g \in \mathcal{G}$ be the set of generators with $|\mathcal{G}| = G$. Variable~$\mathbf{p} \in \mathbb{R}^{G \times 1}$ represents the power production of generators. To define the Power Transfer Distribution Factor~(PTDF)\footnote{The PTDFs measure the change in line power flow induced by injecting one unit of power at any bus and withdrawing it at the slack bus.} matrix, we introduce the line susceptances~$\mathbf{b} \in \mathbb{R}^{L \times 1}$ as well as the line-node incidence matrix~$\mathbf{K} \in \mathbb{Z}^{L \times N}$ and the slack bus adjusted line-node incidence matrix~$\mathbf{K}^\mathrm{sba} \in \mathbb{Z}^{L \times N-1}$, respectively. They define a starting line at the node with \num{1} and an ending line with \num{-1}. We define the slack bus adjusted PTDF matrix of the full system~$\mathbf{PTDF}^\mathrm{sba} \in \mathbb{R}^{L\times N-1}$ by:
    \begin{equation}\label{eq:ptdf}
        \mathbf{PTDF}^\mathrm{sba}=\operatorname{diag}(\mathbf{b})\, \mathbf{K}^\mathrm{sba} \, ((\mathbf{K}^\mathrm{sba})\tran \, \operatorname{diag}(\mathbf{b}) \, \mathbf{K}^\mathrm{sba})^{-1}.
    \end{equation}
    To obtain the PTDF matrix of the full grid~$\mathbf{PTDF} \in \mathbb{R}^{L \times N}$, we add a zero column at the index of the slack bus to the $\mathbf{PTDF}^{sba}$. Additionally, we define the following parameters: generator production costs~$\mathbf{C} \in \mathbb{R}^{G \times 1}$, maximum generation~$\mathbf{\bar P}  \in \mathbb{R}^{G \times 1}$, nodal electricity demand~$\mathbf{D} \in \mathbb{R}^{N \times 1}$, row vectors of ones~$\mathbf{e}^\textrm{G} \in \mathbb{R}^{1 \times G}$ and $\mathbf{e}^\textrm{N} \in \mathbb{R}^{1 \times N}$, maximum line power flow~$\mathbf{\bar T}\in \mathbb{R}^{L \times 1}$, and generator to node mapping~$\boldsymbol{\Gamma}\in \mathbb{Z}^{N \times G}$. The line power flows can be included explicitly as an auxiliary variable in the DC-OPF or calculated ex-post as $\boldsymbol{f} = \mathbf{PTDF}\, \left(\boldsymbol{\Gamma} \, {\mathbf{p}} -  \mathbf{D} \right)$. Then we define the single-period DC-OPF problem by:
    \begin{subequations}\label{eq:DC-OPF}
        \begin{align}
            \min_{\mathbf{p}} \quad & z = \mathbf{C}\tran \, \mathbf{p} & & \label{eq:objective_function}\\
            \textrm{s.t.} \quad & \mathbf{0} \leq \mathbf{p} \leq \mathbf{\bar P}  & &: \boldsymbol{\munderbar{\eta}}, \boldsymbol{\bar \eta}, \label{eq:max_prod} \\
            &  \mathbf{e}^\textrm{N} \, \mathbf{D} = \mathbf{e}^\textrm{G}\,  \mathbf{p} & &:\lambda_{slack}, \label{eq:system_balance_constraint_compact} \\
            & - \mathbf{\bar T} \leq \mathbf{PTDF} \, (\boldsymbol{\Gamma} \, \mathbf{p} - \mathbf{D}) \leq \mathbf{\bar T} & &:\phiLoVec, \phiUpVec, \label{eq:line_limit_constraint_compact}
        \end{align}
    \end{subequations}
    where the objective~\eqref{eq:objective_function} is to minimize the total system operating costs~$z$.
    Constraints~\eqref{eq:max_prod} limit the power production of generators and \eqref{eq:system_balance_constraint_compact} enforce the system power balance. Constraints~\eqref{eq:line_limit_constraint_compact} establish the lower and upper limits of the line power flow. We define the dual variables for Constraints~\eqref{eq:max_prod}--\eqref{eq:line_limit_constraint_compact} after the colon.

\section{Grid Aggregation Methodology}\label{s:appendix_grid_aggregation_methodology}

    Section~\ref{s:nodal_congestion_cost} introduces the distance metrics for grid partitioning that are subject to our comparison. One of these metrics is the novel Network Congestion Price. Secondly, we present the combinations of distance metrics and clustering algorithms in Section~\ref{s:clustering_methods}. Finally, we present a method to evaluate different grid aggregations in Section~\ref{s:evaluation}.

    \subsection{Distance Metrics Based on DC-OPF}\label{s:nodal_congestion_cost}
        Several studies~\cite{singh_reduced_2005, imran_effectiveness_2008, cao_incorporating_2018} propose using the LMPs as a distance metric for grid partitioning. As shown in~\cite{chatzivasileiadis_optimization_2018}, for the model in~\eqref{eq:DC-OPF}, we can derive the LMPs, using the dual variable of the system balance constraint~\eqref{eq:system_balance_constraint_compact} and the dual variables~$\phiLoVec$ and $\phiUpVec$ of the line limit Constraints~\eqref{eq:line_limit_constraint_compact} as:
        \begin{equation}\label{eq:LMP_definition}
            \mathbf{LMP} = \lambda_{slack} + \mathbf{PTDF}\tran \,(\phiUpVec -\phiLoVec),
        \end{equation}
        where $\lambda_{slack}$ denotes the LMP at the slack bus and represents the energy component of the LMPs. The second term represents the influence of line congestions. Consequently, when there is no congestion, the LMPs are the same for all nodes. However, when there is congestion, the term~$(\phiUpVec - \phiLoVec)$ has nonzero entries and the LMPs generally differ at each node. Note, that when multiple lines are congested, the multiplication with the PTDF matrix aggregates the individual impacts captured by the elements of~$(\phiUpVec - \phiLoVec)$. Hence, this aggregation causes a loss of information about the impact of individual line congestions on the LMPs. In fact, considering only the LMPs, it is not even clear how many congestions exist. In Section~\ref{s:results_network_partitioning} we show that the LMPs may be unsuitable as a distance metric for grid aggregation when there is more than one line congested.
        To overcome this problem, we propose using the novel NCPs as a distance metric for grid partitioning for DC-OPF.
        We define the~$\NCP \in \mathbb{R}^{N \times L}$ as:
        \begin{equation}\label{eq:node_cluster_information}
            \NCP = \mathbf{PTDF}\tran \, \operatorname{diag}(\phiUpVec - \phiLoVec).
        \end{equation}
        In contrast to the calculation of the LMPs in Equation~\eqref{eq:LMP_definition}, the $\NCP$  uses a diagonal matrix of the difference of the dual variables~$\phiLoVec$ and~$\phiUpVec$. Consequently, the $\NCP$ does not aggregate the impacts of different line congestions but preserves their individual price contributions.
        Note that while its dimension is generally $N \times L$, all entries of uncongested lines are zero. Therefore, we highlight that the two distance metrics LMP and $\NCP$ lead to identical aggregation results when only a single line is congested. Based on two stylized case studies with multiple line congestions, we will demonstrate in Section~\ref{s:case_study_results} how the higher level of information in the $\NCP$ leads to better results than when using LMPs in the grid partitioning for DC-OPF.

    \subsection{Grid Partitioning Methods}\label{s:clustering_methods}
        In the following, we introduce five different partitioning methods~$k \in \mathcal{K}$ combining the two distance metrics (LMPs, NCPs) with various clustering algorithms. For a user-selected number of nodes $\tilde{N}$ in the aggregated grid, all methods $k$ result in a node-cluster mapping matrix~$\mathbf{M}^\mathrm{nc}_k \in \mathbb{Z}^{\tilde{N} \times N}$, as well as a line mapping matrix~$\mathbf{M}^\mathrm{l}_k \in \mathbb{Z}^{\tilde L \times L}$, assigning nodes and lines of the full grid to those in the aggregated grid. We omit the centroids from the clustering and use only the group assignments to aggregate nodes.
        Note, that we do not reoptimize power flows between iterations of the grid aggregation process. Instead, we always consider the optimal schedules of the full grid when making clustering decisions.
        
        As a benchmark, we use the LMPs as a distance metric in a KMeans algorithm~\cite{imran_effectiveness_2008}. Since this approach does not consider the direct connection of nodes during clustering, the authors in~\cite{cao_incorporating_2018} propose using a Spectral Clustering~(SC) algorithm instead. For completeness, we repeat the mathematical definition of the SC~algorithm in Appendix~\ref{sec:appendix_SC}. However, the SC algorithm fails to capture important line congestions, which we discuss in more detail in Section~\ref{s:results_network_partitioning}. Therefore, we introduce an algorithm based on Agglomerative Clustering~\cite{pedregosa_scikit-learn_nodate}, which updates the distances between iterations such that only adjacent nodes can be grouped, referred to as Adjacent Node Agglomerative Clustering~(ANAC). Appendix~\ref{sec:Appendix_ACA} provides a comprehensive description of the algorithm. In Table~\ref{tab:clustering_methods} we outline the combination of distance metric~(LMP, NCP) and the clustering algorithms~(KMeans, SC, ANAC) used for the comparison of grid partitionings. We use the SC~algorithm only in combination with LMPs because it requires a vector with dimensions $L\times 1$ as an input, which does not match the dimensions of the $\NCP$.
        \begin{table}[!t] \renewcommand{\arraystretch}{1.3} \caption{Combinations of distance metric and clustering algorithm used for grid partitioning indicated by~$\times$.} \label{tab:clustering_methods} 
        \centering 
        \begin{tabular}{l|c|c|c} \hline
                \multirow[b]{2}{*}{\textbf{Distance Metric}} & \multicolumn{3}{c}{\textbf{Clustering Algorithm}}\\\cline{2-4}
                 & KMeans & Spectral Clustering  & ANAC   \\ \hline
                $\mathbf{LMP}$  & $\times$& $\times$ & $\times$  \\
                \hline
                 $\mathbf{NCP}$& $\times$ & - & $\times$    \\
                 \hline
            \end{tabular} 
            \vspace{-0.3cm}
        \end{table}
    \subsection{Evaluation Process}\label{s:evaluation}
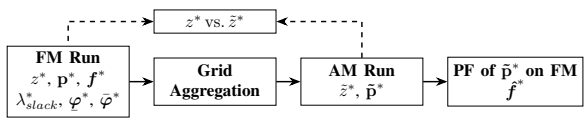
\begin{figure}[b] 
            \centering \resizebox{\columnwidth}{!}{%
    \usetikzlibrary{arrows} 
    \usetikzlibrary{positioning}
    
    \begin{tikzpicture}[
        node distance= 0.5cm and 0cm,
        every path/.style={line width=0.3 mm},
        every node/.style={rectangle, draw=black, fill=white, text=black, align=center, minimum width=2.5cm, line width=0.1mm
        }
    ]
    
        \node (main) {\textbf{FM Run} \\ $z^*$, $\mathbf{p}^*$, $\boldsymbol{f}^*$ \\ $\lambda_{slack}^*$, $\phiLoVec^*$, $\phiUpVec^*$};
        \node (aggregation) [right =of main, xshift=0.5cm] {\textbf{Grid} \\ \textbf{Aggregation}};
        \node (am_run) [right =of aggregation, xshift=0.5cm] {\textbf{AM Run} \\ $\tilde{z}^*$, $\mathbf{\tilde{p}}^*$};
        \node (pf_am) [right =of am_run, xshift=0.5cm] {\textbf{PF of} $\mathbf{\tilde{p}}^*$ \textbf{on FM} \\ $\pfVecOp$};
        \node (comp_ofv)[above = of aggregation] {$z^*\, \mathrm{vs.}\, \tilde{z}^*$};

        \draw [-{Stealth}] (main.east) -- (aggregation.west);
        \draw [-{Stealth}] (main.east)(aggregation.east) -- (am_run.west);
        \draw [-{Stealth}] (main.east)(am_run.east) -- (pf_am.west);

        \draw [dashed, -{Stealth}] (am_run.north) |- (comp_ofv.east);
        \draw [dashed, {Stealth}-] (comp_ofv.west) -| (main.north);
        
    \end{tikzpicture}
}
            \vspace{-0.7cm}
            \caption{Evaluation of grid aggregation methods.}
            \label{fig:evaluation_method}
        \end{figure}

        To evaluate the quality of the grid partitioning methods presented in Table~\ref{tab:clustering_methods}, we use the process illustrated in Figure~\ref{fig:evaluation_method}, which we explain in detail in the following.
        Initially, we solve the DC-OPF problem~\eqref{eq:DC-OPF} for the grid of the full model~(FM) and obtain the optimal solutions for~$z^*, \mathbf{p}^*, \boldsymbol{f}^*, \lambda_{slack}^*, \phiLoVec^*$, and $\phiUpVec^*$. As shown in~\cite{fortenbacher_transmission_2018}, we aggregate the full grid and determine the reduced PTDF matrix according to~\eqref{eq:reduced_PTDF} based on the~$\mathbf{PTDF}$ of the full grid as well as the node mapping~$\mathbf{M}^\mathrm{nc}_k$ and line mapping~$\mathbf{M}^\mathrm{l}_k$ between the full and the reduced grid, respectively.
        \begin{equation}\label{eq:reduced_PTDF}
            \mathbf{PTDF}^{\mathrm{r}} = \mathbf{M}_k^\mathrm{l} \, \mathbf{PTDF}\, {(\mathbf{M}_k^\mathrm{nc}})\tran \, (\mathbf{M}_k^\mathrm{nc}\,{(\mathbf{M}_k^\mathrm{nc}})\tran)^{-1} 
        \end{equation}
        We only assign generators to the nodes in the aggregated grid according to the mapping matrix~$\mathbf{M}^\mathrm{nc}_k$. Then, we solve the DC-OPF model~\eqref{eq:DC-OPF} with the aggregated grid~(AM) and retrieve the optimal solutions for the objective function value~$\tilde{z}^*$ and the unit dispatch~$\tilde{\mathbf{p}}^*$. While the unit dispatch can be compared directly between FM and AM, as there is no reduction of generators, the optimal line power flows cannot be straightforwardly compared due to a reduction in the set of lines. Instead, we use the optimal unit dispatch~$\tilde{\mathbf{p}}^*$ of the AM and the PTDF matrix of the full grid to derive the flows~$\pfVecOp$ that would be attained if~$\tilde{\mathbf{p}}^*$ was the optimal dispatch of the FM:
        \begin{equation}\label{eq:powerFlow_mapping}
            \pfVecOp = \mathbf{PTDF} \left(\boldsymbol{\Gamma} \, \tilde{\mathbf{p}}^* -  \mathbf{D} \right).
        \end{equation}
        Based on the AM and FM results, we define two error metrics for the objective function value and the line limit violations.
        Since the objective function is the sum of all dispatched units multiplied by their cost, we use the objective function as a metric to compare the differences in unit dispatch. Therefore, we define the relative error in the objective function value~$\mathrm{ROVE}(z^*,\tilde{z}^*)$ as:
        \begin{equation}\label{eq:ovf_relative_error}
            \mathrm{ROVE}(z^*,\tilde{z}^*) = \frac{\tilde{z}^* - z^*}{z^*}.
        \end{equation}   
        To evaluate the optimal power flows~$\pfVecOp$ derived from the AM by Equation~\eqref{eq:powerFlow_mapping}, we focus on the violation of line limits and define the maximum relative line limit violation~$\mathrm{MRLLV}(\pfVecOp)$ as:\footnote{In this paper, we always understand the division by a vector element-wise.}
        \begin{equation}\label{eq:powerFlow_relative_error}
            \mathrm{MRLLV}(\pfVecOp) = \mathrm{max}\left( \frac{|\pfVecOp| - \mathbf{\bar T}}{\mathbf{\bar T}}, 0 \right).
        \end{equation}
        In this work, we focus on the quality of the grid aggregation only and postpone an analysis of the solving time of the AMs of different aggregation methods to future research. However, we do compare the time required for the grid partitioning process, referred to as grid partitioning time~(GPT). Which we measure with the results of the FM as a starting point and the resulting node-cluster mapping matrix~$\mathbf{M}^\mathrm{nc}_k$ as the end.

\section{Case Studies \& Results}
\label{s:case_study_results}
    To illustrate the advantages of the proposed $\NCP$ over the LMP distance metric for the case of multiple line congestions, we apply all aggregation methods introduced in Section~\ref{s:clustering_methods} to two stylized case studies and evaluate their performance as explained in Section~\ref{s:evaluation}. In the following Section~\ref{sec:case_study} we introduce the case studies, whereas Sections~\ref{s:results_network_partitioning} and \ref{s:results_grid_aggregation} show the results of the grid partitioning and the performance of the aggregated grid models, respectively.

    \subsection{IEEE RTS 24-Bus System \& IEEE 300-Bus System}
    \label{sec:case_study}
\begin{figure}
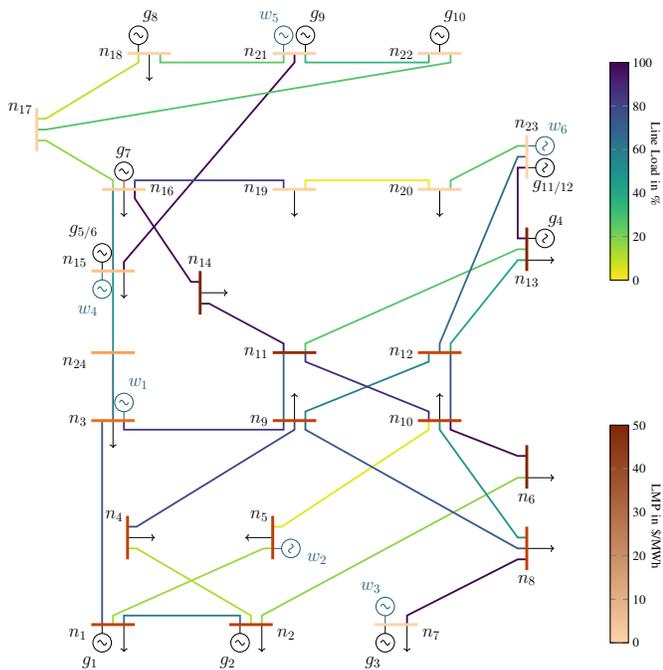

            \centering
                \include{figures/case_study}
                \vspace{-0.7cm}
            \caption{Illustration of the adapted IEEE RTS 24-Bus System~\cite{ordoudis_updated_nodate}. The colors show the optimal line loading $\boldsymbol{\bar f}^*$ and the LMPs.}
            \label{fig:large_case}
            \vspace{-0.3cm}
        \end{figure}

        For an illustrative test case, we use the IEEE RTS 24-Bus System in~\cite{ordoudis_updated_nodate} as a starting point. To create several line congestions in the grid, we increase the installed wind capacity and demand values by a factor of two. For a bigger test case, we use an adapted version of the IEEE 300-Bus System from the \textit{Power Grid Library - OPF}~\cite{babaeinejadsarookolaee_power_2021}, which contains single-period test cases for power flow studies. We increase the demand by a factor of \num{1.2} and add \num{23} wind generators with an installed capacity in the scale of the adjacent line limits to create multiple line congestions. Additionally, we set the minimum production limits of generators to zero as we neglect unit commitment. 
        Furthermore, we perturb the production costs of thermal generators in both case studies so that they are all unique. This modification should foster unique optimization model outcomes, which are desirable for comparing different grid aggregation methods. We further employ backup generators with comparable high costs at all nodes to ensure model feasibility. In this paper, we consider a single timestep only (hour 19 for the IEEE RTS 24-Bus System) and delegate the multi-period case to future research. 
        We implement all models in Python using Pyomo~v6.8.0~\cite{noauthor_pyomo_nodate} and solve them with Gurobi~v11.0.1~\cite{Gurobi2023}. We use Scikit-learn~\cite{pedregosa_scikit-learn_nodate} for clustering. The complete data set and model implementations can be found in~\cite{stockl_grid_2025}. We illustrate the 24-Bus System along with the optimal line loading $\boldsymbol{\bar f}^* = |\boldsymbol{f}^*|/\mathbf{\bar T}$ and the LMPs obtained by~\eqref{eq:LMP_definition} for the FM in Figure~\ref{fig:large_case}.
    \subsection{Results: Grid Partitioning}
    \label{s:results_network_partitioning}
        %
        %
        \begin{figure*}
            \captionsetup[subfigure]{aboveskip=-15pt, margin=1pt}
            \begin{subfigure}{0.19 \textwidth}
            \centering
                \resizebox{0.9\columnwidth}{!}{%

     \begin{circuitikz}  
       
            \fill [IEE_blue!40!white] (2.5,4) -- (2.5,7.5) -- (4.9, 9.036) -- (4.9, 10.9625) -- (11, 10.9625) -- (11, 13) -- (14.3, 13) -- (14.3, 7.75) -- (15.5, 7.75)  -- (15.5,5.5) -- (11, 5.5) -- (8.5,4) -- (2.5,4);
            
            \fill [IEE_red!40!white]  (2.5,13.225) -- (4.7, 13.225) -- (4.7, 14.5) -- (6, 16.02) -- (15.5, 16.02) -- (15.5, 20.75) -- (2.5, 20.75) -- (2.5, 20) -- (1, 19.15) -- (1, 17.15) -- (2.5, 16) -- (2.5, 13.225) ;
            
            \fill [IEE_green!40!white]  (2.5, 13.025) -- (4.7, 13.025) -- (4.7, 11.1625) -- (2.5, 11.1625)  ;
            \fill [IEE_orange!40!white] (2.5, 10.9625) -- (4.7, 10.9625) -- (4.7, 9.236) -- (2.5, 7.75) -- (2.5, 10.9625)  ;
            \fill [IEE_red!40!white] (8.85, 4) -- (11.05, 5.3) -- (15.5, 5.3) -- (15.5, 4) --  (8.85, 4)  ;
            
            \fill [IEE_light_blue!40!white] (6.13, 15.82) -- (15.5, 15.82) -- (15.5, 7.95) -- (14.5, 7.95) -- (14.5, 13.2) -- (10.8, 13.2) --(10.8, 11.1625) -- (4.9, 11.1625) -- (4.9, 14.43);

            \node[draw, circle, inner sep=2pt] at (14.5 ,19.5) {\Huge 1};
            \node[draw, circle, inner sep=2pt] at (3    , 12.55) {\Huge 2};
            \node[draw, circle, inner sep=2pt] at (3.85 , 9.3) {\Huge 3};
            \node[draw, circle, inner sep=2pt] at (13.5 , 11) {\Huge 4};
            \node[draw, circle, inner sep=2pt] at (9 , 14) {\Huge 5};
            \node[draw, circle, inner sep=2pt] at (14.5 ,4.65) {\Huge 1};

            \draw (3+0.9, 4.5)      -- +(0, 0.25)   -- (8-0.9, 4.5+0.25)        -- +(0, -0.25);    
            \draw (3+0.3, 10.125)   -- +(0, -0.25)  -- (3+0.3, 4.5+0.25)        -- +(0, -0.25);     
            \draw (3+0.6, 4.5)      -- +(0, 0.25)   -- (8-0.25, 7.5-0.9)        -- +(0.25, 0);      
            \draw (8-0.6, 4.5)      -- +(0, 0.25)   -- (4.25, 7.5-0.9)          -- +(-0.25, 0);     
            \draw (8-0.3, 4.5)      -- +(0, 0.25)   -- (15-0.25, 8.25+0.3)      -- +(0.25, 0);   
            \draw (3+0.9, 10.125)   -- +(0, -0.25)  -- (8+0.3, 10.125-0.25)     -- +(0, 0.25);     
            \draw (4-0.4, 10.125)                   -- (3+0.6, 12) [line width=0.5mm];                                 
            \draw (4, 7.5-0.3)      -- +(0.25, 0)   -- (8+0.6, 10.125-0.25)     -- +(0, 0.25);       
            \draw (8, 7.5-0.3)      -- +(0.25, 0)   -- (12+0.3, 10.125-0.25)    -- +(0, 0.25);    
            \draw [color=red] (15, 8.25+0.9)    -- +(-0.25, 0)  -- (12+0.9, 10.125-0.25)    -- +(0, 0.25); 
            \draw [color=red](12-0.3, 4.5)     -- +(0, 0.25)   -- (15-0.25, 6+0.3)         -- +(0.25, 0); 
            \draw (15, 6+0.6)       -- +(-0.25, 0)  -- (8+0.9, 10.125-0.25)     -- +(0, 0.25);   
            \draw (15, 6+0.9)       -- +(-0.25, 0)  -- (12+0.6, 10.125-0.25)    -- +(0, 0.25); 
            \draw (8+0.3, 10.125)                                               -- (8+0.3, 12);  
            \draw (8+0.9, 10.125)   -- +(0, 0.25)   -- (12+0.3, 12-0.25)        -- +(0, 0.25);      
            \draw (12+0.3, 10.125)  -- +(0, 0.25)   -- (8+0.9, 12-0.25)         -- +(0, 0.25);      
            \draw (12+0.9, 10.125)                  -- (12+0.9, 12);                                
            \draw (8+0.9, 12)       -- +(0, 0.25)   -- (15-0.25, 14.25+0.6)     -- +(0.25, 0);      
            \draw (6, 14.25-0.9)    -- +(0.25, 0)   -- (8+0.3, 12+0.25)         -- +(0, -0.25);     
            \draw (12+0.9, 12)      -- +(0, 0.25)   -- (15-0.25, 14.25+0.3)     -- +(0.25, 0);      
            \draw (12+0.6, 12)      -- +(0, 0.25)   -- (15-0.25, 18-0.6)        -- +(0.25, 0);      
            \draw [color=red](15, 14.25+0.9)   -- +(-0.25, 0)  -- (15-0.25, 18-0.9)        -- +(0.25, 0);      
            \draw [color=red] (6, 14.25-0.3)    -- +(-0.25, 0)  -- (4.5-0.3, 16.5-0.25)     -- +(0, 0.25);      
            \draw (3+0.6, 14.25)    -- +(0, 0.25)   -- (4.5-0.9, 16.5-0.25)     -- +(0, 0.25); 
            \draw [color=red] (3+0.9, 14.25)    -- +(0, 0.25)   -- (8+0.6, 20.25-0.25)      -- +(0, 0.25); 
            \draw (3+0.6, 14.25)                    -- (3+0.6, 12);                            
            \draw (1.5, 18.75-0.9)  -- +(0.25, 0)   -- (4.5-0.9, 16.5+0.25)     -- +(0, -0.25);  
            \draw (4.5-0.3, 16.5)   -- +(0, 0.25)   -- (8+0.3, 16.5+0.25)       -- +(0, -0.25);  
            \draw (1.5, 18.75-0.3)  -- +(0.25, 0)   -- (4+0.3, 20.25-0.25)      -- +(0, 0.25);  
            \draw (1.5, 18.75-0.6)  -- +(0.25, 0)   -- (12+0.9, 20.25-0.25)     -- +(0, 0.25);  
            \draw (4+0.9, 20.25)    -- +(0, -0.25)  -- (8+0.3, 20.25-0.25)      -- +(0, 0.25);  
            \draw (8+0.9, 16.5)     -- +(0, 0.25)   -- (12+0.3, 16.5+0.25)      -- +(0, -0.25); 
            \draw (12+0.9, 16.5)    -- +(0, 0.25)   -- (15-0.25, 18-0.3)        -- +(0.25, 0);  
            \draw (8+0.9, 20.25)    -- +(0, -0.25)  -- (12+0.3, 20.25-0.25)     -- +(0, 0.25); 

                
            \foreach \start/\name/\direction/\rot/\shift  in {
                {3, 4.5}/n_{1}/north east/0/{1.2,0},
                {8, 4.5}/n_{2}/north west/0/{-1.2,0},
                {3, 10.125}/n_{3}/east/0/{1.2,0},       
                {4, 7.5}/n_{4}/east/0/{0,-1.2},      
                {8, 7.5}/n_{5}/east/0/{0,-1.2},       
                {15, 8.25}/n_{6}/north/0/{0,1.2},  
                {12, 4.5}/n_{7}/north west/0/{-1.2,0},
                {15, 6}/n_{8}/north/0/{0,1.2},    
                {8, 10.125}/n_{9}/east/0/{1.2,0},        
                {12, 10.125}/n_{10}/east/0/{1.2,0},     
                {8, 12}/n_{11}/east/0/{1.2,0},
                {12, 12}/n_{12}/east/0/{1.2,0},
                {15, 14.25}/n_{13}/north/0/{0,1.2},
                {6, 14.25}/n_{14}/south/0/{0,-1.2},
                {3, 14.25}/n_{15}/south east/0/{1.2,0},
                {4.5, 16.5}/n_{16}/west/0/{-1.2,0}, 
                {1.5, 18.75}/n_{17}/east/0/{0,-1.2},
                {4, 20.25}/n_{18}/east/0/{1.2,0},
                {8, 16.5}/n_{19}/east/0/{1.2,0}, 
                {12, 16.5}/n_{20}/east/0/{1.2,0}, 
                {8, 20.25}/n_{21}/east/0/{1.2,0},
                {12, 20.25}/n_{22}/east/0/{1.2,0},
                {15, 18}/n_{23}/south/0/{0,-1.2},
                {3, 12}/n_{24}/north east/0/{1.2,0}}
            {
                \node[anchor=\direction] (\name) at (\start) {};
                \draw [ultra thick] (\start) -- ($(\start)+(\shift)$);
            }

    \end{circuitikz}
    
    
} 
                \caption{LMP-KMeans}
                \label{fig:NP_lmps_kmeans}
            \end{subfigure}
            \begin{subfigure}{0.19 \textwidth}
            \centering
                \resizebox{0.9\columnwidth}{!}{%

     \begin{circuitikz}  
        

            \fill [IEE_red!40!white]  (17.5, 20.75)  -- (30.5, 20.75) -- (30.5, 18.6)  -- (17.5, 17.15) -- (16, 17.15) -- (16, 19.15) -- (17.5, 20);

           \fill [IEE_blue!40!white] (17.5,4) -- (17.5,7.5) -- (19.9, 9.036) -- (19.9, 10.9625) -- (28, 10.9625) -- (30.5, 10.9625) -- (30.5, 7.95) -- (23.5,4) -- (19.5,4);

            

            \fill [IEE_green!40!white] (17.5, 16.95) -- (30.5, 18.4) -- (30.5, 11.1625) -- (19.9, 11.1625) -- (19.9, 15.475) -- (17.5, 15.475) ;
            
            \fill [IEE_orange!40!white] (17.5, 15.275) -- (19.7, 15.275) -- (19.7, 9.236) -- (17.5, 7.75)  ;
            
            \fill [IEE_light_blue!40!white] (23.85, 4) -- (26.05, 5.3) -- (30.5, 7.75) -- (30.5, 4) --  (23.85, 4)  ;


            \node[draw, circle, inner sep=2pt] at (29.5 ,19.5) {\Huge 1};
            \node[draw, circle, inner sep=2pt] at (25   , 14.55) {\Huge 2};
            \node[draw, circle, inner sep=2pt] at (18.85 , 9.3) {\Huge 3};
            \node[draw, circle, inner sep=2pt] at (29.5 , 10) {\Huge 4};
            \node[draw, circle, inner sep=2pt] at (29.5 ,4.65) {\Huge 5};

            \draw (3   + 15 + 0.9, 4.5      )  -- +(0, 0.25)   -- (8    + 15 - 0.9 , 4.5+0.25)        -- +(0, -0.25);    
            \draw (3   + 15 + 0.3, 10.125   )  -- +(0, -0.25)  -- (3    + 15 + 0.3 , 4.5+0.25)        -- +(0, -0.25);     
            \draw (8   + 15 - 0.6, 4.5      )  -- +(0, 0.25)   -- (4.25 + 15       , 7.5-0.9)          -- +(-0.25, 0);     
            \draw (3   + 15 + 0.6, 4.5      )  -- +(0, 0.25)   -- (8    + 15 - 0.25, 7.5-0.9)        -- +(0.25, 0);      
            \draw (8   + 15 - 0.3, 4.5      )  -- +(0, 0.25)   -- (15   + 15 - 0.25, 8.25+0.3)      -- +(0.25, 0);   
            \draw (3   + 15 + 0.9, 10.125   )  -- +(0, -0.25)  -- (8    + 15 + 0.3 , 10.125-0.25)     -- +(0, 0.25);     
            \draw (4   + 15 - 0.4, 10.125   )                  -- (3    + 15 + 0.6 , 12);                                 
            \draw (4   + 15      , 7.5-0.3  )  -- +(0.25, 0)   -- (8    + 15 + 0.6 , 10.125-0.25)     -- +(0, 0.25);       
            \draw (8   + 15      , 7.5-0.3  )  -- +(0.25, 0)   -- (12   + 15 + 0.3 , 10.125-0.25)    -- +(0, 0.25);    
            \draw [color=red] (15  + 15      , 8.25+0.9 )  -- +(-0.25, 0)  -- (12   + 15 + 0.9 , 10.125-0.25)    -- +(0, 0.25); 
            \draw [color=red](12  + 15 - 0.3, 4.5      )  -- +(0, 0.25)   -- (15   + 15 - 0.25, 6+0.3)         -- +(0.25, 0); 
            \draw (15  + 15      , 6+0.6    )  -- +(-0.25, 0)  -- (8    + 15 + 0.9 , 10.125-0.25)     -- +(0, 0.25);   
            \draw (15  + 15      , 6+0.9    )  -- +(-0.25, 0)  -- (12   + 15 + 0.6 , 10.125-0.25)    -- +(0, 0.25); 
            \draw (8   + 15 + 0.3, 10.125   )                  -- (8    + 15 +0.3, 12);  
            \draw (8   + 15 + 0.9, 10.125   )  -- +(0, 0.25)   -- (12   + 15 + 0.3 , 12-0.25)        -- +(0, 0.25);      
            \draw (12  + 15 + 0.3, 10.125   )  -- +(0, 0.25)   -- (8    + 15 + 0.9 , 12-0.25)         -- +(0, 0.25);      
            \draw (12  + 15 + 0.9, 10.125   )                  -- (12   + 15 + 0.9 , 12);                                
            \draw (8   + 15 + 0.9, 12       )  -- +(0, 0.25)   -- (15   + 15 - 0.25, 14.25+0.6)     -- +(0.25, 0);      
            \draw (6   + 15      , 14.25-0.9)  -- +(0.25, 0)   -- (8    + 15 + 0.3 , 12+0.25)         -- +(0, -0.25);     
            \draw (12  + 15 + 0.9, 12       )  -- +(0, 0.25)   -- (15   + 15 - 0.25, 14.25+0.3)     -- +(0.25, 0);      
            \draw (12  + 15 + 0.6, 12       )  -- +(0, 0.25)   -- (15   + 15 - 0.25, 18-0.6)        -- +(0.25, 0);      
            \draw [color=red](15  + 15      , 14.25+0.9)  -- +(-0.25, 0)  -- (15   + 15 - 0.25, 18-0.9)        -- +(0.25, 0);      
            \draw [color=red](6   + 15      , 14.25-0.3)  -- +(-0.25, 0)  -- (4.5  + 15 - 0.3 , 16.5-0.25)     -- +(0, 0.25);      
            \draw (3   + 15 + 0.6, 14.25    )  -- +(0, 0.25)   -- (4.5  + 15 - 0.9 , 16.5-0.25)     -- +(0, 0.25); 
            \draw [color=red] (3   + 15 + 0.9, 14.25    )  -- +(0, 0.25)   -- (8    + 15 + 0.6 , 20.25-0.25)      -- +(0, 0.25); 
            \draw (3   + 15 + 0.6, 14.25    )                  -- (3    + 15 + 0.6 , 12);                            
            \draw (1.5 + 15      , 18.75-0.9)  -- +(0.25, 0)   -- (4.5  + 15 - 0.9 , 16.5+0.25)     -- +(0, -0.25);  
            \draw (4.5 + 15 - 0.3, 16.5     )  -- +(0, 0.25)   -- (8    + 15 + 0.3 , 16.5+0.25)       -- +(0, -0.25);  
            \draw (1.5 + 15      , 18.75-0.3)  -- +(0.25, 0)   -- (4    + 15 + 0.3 , 20.25-0.25)      -- +(0, 0.25);  
            \draw (1.5 + 15      , 18.75-0.6)  -- +(0.25, 0)   -- (12   + 15 + 0.9 , 20.25-0.25)     -- +(0, 0.25);  
            \draw (4   + 15 + 0.9, 20.25    )  -- +(0, -0.25)  -- (8    + 15 + 0.3 , 20.25-0.25)      -- +(0, 0.25);  
            \draw (8   + 15 + 0.9, 16.5     )  -- +(0, 0.25)   -- (12   + 15 + 0.3 , 16.5+0.25)      -- +(0, -0.25); 
            \draw (12  + 15 + 0.9, 16.5     )  -- +(0, 0.25)   -- (15   + 15 - 0.25, 18-0.3)        -- +(0.25, 0);  
            \draw (8   + 15 + 0.9, 20.25    )  -- +(0, -0.25)  -- (12   + 15 + 0.3 , 20.25-0.25)     -- +(0, 0.25); 

                
            \foreach \start/\name/\direction/\rot/\shift  in {
                {  3 + 15 , 4.5}/n_{1}/north east/0/{1.2,0},
                {  8 + 15 , 4.5}/n_{2}/north west/0/{-1.2,0},
                {  3 + 15 , 10.125}/n_{3}/east/0/{1.2,0},       
                {  4 + 15 , 7.5}/n_{4}/east/0/{0,-1.2},      
                {  8 + 15 , 7.5}/n_{5}/east/0/{0,-1.2},       
                { 15 + 15 , 8.25}/n_{6}/north/0/{0,1.2},  
                { 12 + 15 , 4.5}/n_{7}/north west/0/{-1.2,0},
                { 15 + 15 , 6}/n_{8}/north/0/{0,1.2},    
                {  8 + 15 , 10.125}/n_{9}/east/0/{1.2,0},        
                { 12 + 15 , 10.125}/n_{10}/east/0/{1.2,0},     
                {8   + 15 , 12}/n_{11}/east/0/{1.2,0},
                {12  + 15 , 12}/n_{12}/east/0/{1.2,0},
                {15  + 15 , 14.25}/n_{13}/north/0/{0,1.2},
                {6   + 15 , 14.25}/n_{14}/south/0/{0,-1.2},
                {3   + 15 , 14.25}/n_{15}/south east/0/{1.2,0},
                {4.5 + 15 , 16.5}/n_{16}/west/0/{-1.2,0}, 
                {1.5 + 15 , 18.75}/n_{17}/east/0/{0,-1.2},
                {4   + 15 , 20.25}/n_{18}/east/0/{1.2,0},
                {8   + 15 , 16.5}/n_{19}/east/0/{1.2,0}, 
                {12  + 15 , 16.5}/n_{20}/east/0/{1.2,0}, 
                {8   + 15 , 20.25}/n_{21}/east/0/{1.2,0},
                {12  + 15 , 20.25}/n_{22}/east/0/{1.2,0},
                {15  + 15 , 18}/n_{23}/south/0/{0,-1.2},
                {3   + 15 , 12}/n_{24}/north east/0/{1.2,0}}
            {
                \node[anchor=\direction] (\name) at (\start) {};
                \draw [ultra thick] (\start) -- ($(\start)+(\shift)$);
            }
        

    \end{circuitikz}
    
    
} 
                \caption{LMP-SC}
                \label{fig:NP_lmps_spectral}
            \end{subfigure}
            \begin{subfigure}{0.19 \textwidth}
            \centering
                \resizebox{0.9\columnwidth}{!}{%

     \begin{circuitikz}  

            \fill [IEE_blue!40!white] (17.5,4) -- (17.5,7.5) -- (19.9, 9.036) -- (19.9, 14.4) -- (21.15, 15.82) -- (30.5,15.82)  -- (30.5,5.5) -- (26, 5.5) -- (23.5,4) -- (17.5,4);
            \fill [IEE_red!40!white]  (17.5,13.225) -- (19.7, 13.225) -- (19.7, 14.5) -- (21, 16.02) -- (30.5, 16.02) -- (30.5, 20.75) -- (17.5, 20.75) -- (17.5, 20) -- (16, 19.15) -- (16, 17.15) -- (17.5, 16) -- (17.5, 13.225) ;
            \fill [IEE_green!40!white]  (17.5, 13.025) -- (19.7, 13.025) -- (19.7, 11.1625) -- (17.5, 11.1625) -- (17.5, 13.025)  ;
            \fill [IEE_orange!40!white] (17.5, 10.9625) -- (19.7, 10.9625) -- (19.7, 9.236) -- (17.5, 7.75) -- (17.5, 10.9625)  ;
            \fill [IEE_light_blue!40!white] (23.85, 4) -- (26.05, 5.3) -- (30.5, 5.3) -- (30.5, 4) --  (23.85, 4)  ;

            \node[draw, circle, inner sep=2pt] at (14.5 + 15 ,19.5) {\Huge 1};
            \node[draw, circle, inner sep=2pt] at (3   + 15 , 12.55) {\Huge 2};
            \node[draw, circle, inner sep=2pt] at (3.85  + 15 , 9.3) {\Huge 3};
            \node[draw, circle, inner sep=2pt] at (14.5 + 15 ,12) {\Huge 4};
            \node[draw, circle, inner sep=2pt] at (14.5 + 15 ,4.65) {\Huge 5};

            \draw (3   + 15 + 0.9, 4.5      )  -- +(0, 0.25)   -- (8    + 15 - 0.9 , 4.5+0.25)        -- +(0, -0.25);    
            \draw (3   + 15 + 0.3, 10.125   )  -- +(0, -0.25)  -- (3    + 15 + 0.3 , 4.5+0.25)        -- +(0, -0.25);     
            \draw (8   + 15 - 0.6, 4.5      )  -- +(0, 0.25)   -- (4.25 + 15       , 7.5-0.9)          -- +(-0.25, 0);     
            \draw (3   + 15 + 0.6, 4.5      )  -- +(0, 0.25)   -- (8    + 15 - 0.25, 7.5-0.9)        -- +(0.25, 0);      
            \draw (8   + 15 - 0.3, 4.5      )  -- +(0, 0.25)   -- (15   + 15 - 0.25, 8.25+0.3)      -- +(0.25, 0);   
            \draw (3   + 15 + 0.9, 10.125   )  -- +(0, -0.25)  -- (8    + 15 + 0.3 , 10.125-0.25)     -- +(0, 0.25);     
            \draw (4   + 15 - 0.4, 10.125   )                  -- (3    + 15 + 0.6 , 12);                                 
            \draw (4   + 15      , 7.5-0.3  )  -- +(0.25, 0)   -- (8    + 15 + 0.6 , 10.125-0.25)     -- +(0, 0.25);       
            \draw (8   + 15      , 7.5-0.3  )  -- +(0.25, 0)   -- (12   + 15 + 0.3 , 10.125-0.25)    -- +(0, 0.25);    
            \draw [color=red] (15  + 15      , 8.25+0.9 )  -- +(-0.25, 0)  -- (12   + 15 + 0.9 , 10.125-0.25)    -- +(0, 0.25); 
            \draw [color=red](12  + 15 - 0.3, 4.5      )  -- +(0, 0.25)   -- (15   + 15 - 0.25, 6+0.3)         -- +(0.25, 0); 
            \draw (15  + 15      , 6+0.6    )  -- +(-0.25, 0)  -- (8    + 15 + 0.9 , 10.125-0.25)     -- +(0, 0.25);   
            \draw (15  + 15      , 6+0.9    )  -- +(-0.25, 0)  -- (12   + 15 + 0.6 , 10.125-0.25)    -- +(0, 0.25); 
            \draw (8   + 15 + 0.3, 10.125   )                  -- (8    + 15 +0.3, 12);  
            \draw (8   + 15 + 0.9, 10.125   )  -- +(0, 0.25)   -- (12   + 15 + 0.3 , 12-0.25)        -- +(0, 0.25);      
            \draw (12  + 15 + 0.3, 10.125   )  -- +(0, 0.25)   -- (8    + 15 + 0.9 , 12-0.25)         -- +(0, 0.25);      
            \draw (12  + 15 + 0.9, 10.125   )                  -- (12   + 15 + 0.9 , 12);                                
            \draw (8   + 15 + 0.9, 12       )  -- +(0, 0.25)   -- (15   + 15 - 0.25, 14.25+0.6)     -- +(0.25, 0);      
            \draw (6   + 15      , 14.25-0.9)  -- +(0.25, 0)   -- (8    + 15 + 0.3 , 12+0.25)         -- +(0, -0.25);     
            \draw (12  + 15 + 0.9, 12       )  -- +(0, 0.25)   -- (15   + 15 - 0.25, 14.25+0.3)     -- +(0.25, 0);      
            \draw (12  + 15 + 0.6, 12       )  -- +(0, 0.25)   -- (15   + 15 - 0.25, 18-0.6)        -- +(0.25, 0);      
            \draw [color=red](15  + 15      , 14.25+0.9)  -- +(-0.25, 0)  -- (15   + 15 - 0.25, 18-0.9)        -- +(0.25, 0);      
            \draw [color=red](6   + 15      , 14.25-0.3)  -- +(-0.25, 0)  -- (4.5  + 15 - 0.3 , 16.5-0.25)     -- +(0, 0.25);      
            \draw (3   + 15 + 0.6, 14.25    )  -- +(0, 0.25)   -- (4.5  + 15 - 0.9 , 16.5-0.25)     -- +(0, 0.25); 
            \draw [color=red] (3   + 15 + 0.9, 14.25    )  -- +(0, 0.25)   -- (8    + 15 + 0.6 , 20.25-0.25)      -- +(0, 0.25); 
            \draw (3   + 15 + 0.6, 14.25    )                  -- (3    + 15 + 0.6 , 12);                            
            \draw (1.5 + 15      , 18.75-0.9)  -- +(0.25, 0)   -- (4.5  + 15 - 0.9 , 16.5+0.25)     -- +(0, -0.25);  
            \draw (4.5 + 15 - 0.3, 16.5     )  -- +(0, 0.25)   -- (8    + 15 + 0.3 , 16.5+0.25)       -- +(0, -0.25);  
            \draw (1.5 + 15      , 18.75-0.3)  -- +(0.25, 0)   -- (4    + 15 + 0.3 , 20.25-0.25)      -- +(0, 0.25);  
            \draw (1.5 + 15      , 18.75-0.6)  -- +(0.25, 0)   -- (12   + 15 + 0.9 , 20.25-0.25)     -- +(0, 0.25);  
            \draw (4   + 15 + 0.9, 20.25    )  -- +(0, -0.25)  -- (8    + 15 + 0.3 , 20.25-0.25)      -- +(0, 0.25);  
            \draw (8   + 15 + 0.9, 16.5     )  -- +(0, 0.25)   -- (12   + 15 + 0.3 , 16.5+0.25)      -- +(0, -0.25); 
            \draw (12  + 15 + 0.9, 16.5     )  -- +(0, 0.25)   -- (15   + 15 - 0.25, 18-0.3)        -- +(0.25, 0);  
            \draw (8   + 15 + 0.9, 20.25    )  -- +(0, -0.25)  -- (12   + 15 + 0.3 , 20.25-0.25)     -- +(0, 0.25); 

                
            \foreach \start/\name/\direction/\rot/\shift  in {
                {  3 + 15 , 4.5}/n_{1}/north east/0/{1.2,0},
                {  8 + 15 , 4.5}/n_{2}/north west/0/{-1.2,0},
                {  3 + 15 , 10.125}/n_{3}/east/0/{1.2,0},       
                {  4 + 15 , 7.5}/n_{4}/east/0/{0,-1.2},      
                {  8 + 15 , 7.5}/n_{5}/east/0/{0,-1.2},       
                { 15 + 15 , 8.25}/n_{6}/north/0/{0,1.2},  
                { 12 + 15 , 4.5}/n_{7}/north west/0/{-1.2,0},
                { 15 + 15 , 6}/n_{8}/north/0/{0,1.2},    
                {  8 + 15 , 10.125}/n_{9}/east/0/{1.2,0},        
                { 12 + 15 , 10.125}/n_{10}/east/0/{1.2,0},     
                {8   + 15 , 12}/n_{11}/east/0/{1.2,0},
                {12  + 15 , 12}/n_{12}/east/0/{1.2,0},
                {15  + 15 , 14.25}/n_{13}/north/0/{0,1.2},
                {6   + 15 , 14.25}/n_{14}/south/0/{0,-1.2},
                {3   + 15 , 14.25}/n_{15}/south east/0/{1.2,0},
                {4.5 + 15 , 16.5}/n_{16}/west/0/{-1.2,0}, 
                {1.5 + 15 , 18.75}/n_{17}/east/0/{0,-1.2},
                {4   + 15 , 20.25}/n_{18}/east/0/{1.2,0},
                {8   + 15 , 16.5}/n_{19}/east/0/{1.2,0}, 
                {12  + 15 , 16.5}/n_{20}/east/0/{1.2,0}, 
                {8   + 15 , 20.25}/n_{21}/east/0/{1.2,0},
                {12  + 15 , 20.25}/n_{22}/east/0/{1.2,0},
                {15  + 15 , 18}/n_{23}/south/0/{0,-1.2},
                {3   + 15 , 12}/n_{24}/north east/0/{1.2,0}}
            {
                \node[anchor=\direction] (\name) at (\start) {};
                \draw [ultra thick] (\start) -- ($(\start)+(\shift)$);
            }
        

    \end{circuitikz}
    
    
} 
                \caption{LMP-ANAC}
                \label{fig:NP_lmps_AHC}
            \end{subfigure}
            \begin{subfigure}{0.19 \textwidth}
            \centering
                \resizebox{.9\columnwidth}{!}{%

     \begin{circuitikz}  

            \fill [IEE_blue!40!white] (2.5,4) -- (2.5,10.5) -- (15.5,15.82)  -- (15.5,5.5) -- (11, 5.5) -- (8.5,4) -- (2.5,4);
            \fill [IEE_red!40!white]  (2.5,10.7) -- (4.7, 11.6) -- (4.7, 14.5) -- (8.3, 18.6) -- (15.5, 18.6) -- (15.5, 20.75) -- (2.5, 20.75) -- (2.5, 20) -- (1, 19.15) -- (1, 17.15) -- (2.5, 16) -- (2.5, 10.7) ;
            \fill [IEE_green!40!white]  (15.5,16.02) -- (6.3, 16.02) -- (8.4, 18.4) -- (15.5, 18.4) -- (15.5, 16.02)  ;
            \fill [IEE_orange!40!white] (4.9, 11.681) -- (4.9, 14.4) -- (6.15, 15.82) -- (7, 15.82) -- (7, 12.54) -- (4.9, 11.681)   ;
            \fill [IEE_light_blue!40!white] (8.85, 4) -- (11.05, 5.3) -- (15.5, 5.3) -- (15.5, 4) --  (8.85, 4)  ;

            \node[draw, circle, inner sep=2pt] at (14.5,19.5) {\Huge 1};
            \node[draw, circle, inner sep=2pt] at (12, 17.7) {\Huge 2};
            \node[draw, circle, inner sep=2pt] at (6.3,15) {\Huge 3};
            \node[draw, circle, inner sep=2pt] at (14.5,12) {\Huge 4};
            \node[draw, circle, inner sep=2pt] at (14.5,4.65) {\Huge 5};

            
            \draw (3+0.9, 4.5)      -- +(0, 0.25)   -- (8-0.9, 4.5+0.25)        -- +(0, -0.25);    
            \draw (3+0.3, 10.125)   -- +(0, -0.25)  -- (3+0.3, 4.5+0.25)        -- +(0, -0.25);     
            \draw (3+0.6, 4.5)      -- +(0, 0.25)   -- (8-0.25, 7.5-0.9)        -- +(0.25, 0);      
            \draw (8-0.6, 4.5)      -- +(0, 0.25)   -- (4.25, 7.5-0.9)          -- +(-0.25, 0);     
            \draw (8-0.3, 4.5)      -- +(0, 0.25)   -- (15-0.25, 8.25+0.3)      -- +(0.25, 0);   
            \draw (3+0.9, 10.125)   -- +(0, -0.25)  -- (8+0.3, 10.125-0.25)     -- +(0, 0.25);     
            \draw (4-0.4, 10.125)                   -- (3+0.6, 12);                                 
            \draw (4, 7.5-0.3)      -- +(0.25, 0)   -- (8+0.6, 10.125-0.25)     -- +(0, 0.25);       
            \draw (8, 7.5-0.3)      -- +(0.25, 0)   -- (12+0.3, 10.125-0.25)    -- +(0, 0.25);    
            \draw [color=red] (15, 8.25+0.9)    -- +(-0.25, 0)  -- (12+0.9, 10.125-0.25)    -- +(0, 0.25); 
            \draw [color=red](12-0.3, 4.5)     -- +(0, 0.25)   -- (15-0.25, 6+0.3)         -- +(0.25, 0); 
            \draw (15, 6+0.6)       -- +(-0.25, 0)  -- (8+0.9, 10.125-0.25)     -- +(0, 0.25);   
            \draw (15, 6+0.9)       -- +(-0.25, 0)  -- (12+0.6, 10.125-0.25)    -- +(0, 0.25); 
            \draw (8+0.3, 10.125)                                               -- (8+0.3, 12);  
            \draw (8+0.9, 10.125)   -- +(0, 0.25)   -- (12+0.3, 12-0.25)        -- +(0, 0.25);      
            \draw (12+0.3, 10.125)  -- +(0, 0.25)   -- (8+0.9, 12-0.25)         -- +(0, 0.25);      
            \draw (12+0.9, 10.125)                  -- (12+0.9, 12);                                
            \draw (8+0.9, 12)       -- +(0, 0.25)   -- (15-0.25, 14.25+0.6)     -- +(0.25, 0);      
            \draw (6, 14.25-0.9)    -- +(0.25, 0)   -- (8+0.3, 12+0.25)         -- +(0, -0.25);     
            \draw (12+0.9, 12)      -- +(0, 0.25)   -- (15-0.25, 14.25+0.3)     -- +(0.25, 0);      
            \draw (12+0.6, 12)      -- +(0, 0.25)   -- (15-0.25, 18-0.6)        -- +(0.25, 0);      
            \draw [color=red](15, 14.25+0.9)   -- +(-0.25, 0)  -- (15-0.25, 18-0.9)        -- +(0.25, 0);      
            \draw [color=red] (6, 14.25-0.3)    -- +(-0.25, 0)  -- (4.5-0.3, 16.5-0.25)     -- +(0, 0.25);      
            \draw (3+0.6, 14.25)    -- +(0, 0.25)   -- (4.5-0.9, 16.5-0.25)     -- +(0, 0.25); 
            \draw [color=red] (3+0.9, 14.25)    -- +(0, 0.25)   -- (8+0.6, 20.25-0.25)      -- +(0, 0.25); 
            \draw (3+0.6, 14.25)                    -- (3+0.6, 12);                            
            \draw (1.5, 18.75-0.9)  -- +(0.25, 0)   -- (4.5-0.9, 16.5+0.25)     -- +(0, -0.25);  
            \draw (4.5-0.3, 16.5)   -- +(0, 0.25)   -- (8+0.3, 16.5+0.25)       -- +(0, -0.25);  
            \draw (1.5, 18.75-0.3)  -- +(0.25, 0)   -- (4+0.3, 20.25-0.25)      -- +(0, 0.25);  
            \draw (1.5, 18.75-0.6)  -- +(0.25, 0)   -- (12+0.9, 20.25-0.25)     -- +(0, 0.25);  
            \draw (4+0.9, 20.25)    -- +(0, -0.25)  -- (8+0.3, 20.25-0.25)      -- +(0, 0.25);  
            \draw (8+0.9, 16.5)     -- +(0, 0.25)   -- (12+0.3, 16.5+0.25)      -- +(0, -0.25); 
            \draw (12+0.9, 16.5)    -- +(0, 0.25)   -- (15-0.25, 18-0.3)        -- +(0.25, 0);  
            \draw (8+0.9, 20.25)    -- +(0, -0.25)  -- (12+0.3, 20.25-0.25)     -- +(0, 0.25); 

                
            \foreach \start/\name/\direction/\rot/\shift  in {
                {3, 4.5}/n_{1}/north east/0/{1.2,0},
                {8, 4.5}/n_{2}/north west/0/{-1.2,0},
                {3, 10.125}/n_{3}/east/0/{1.2,0},       
                {4, 7.5}/n_{4}/east/0/{0,-1.2},      
                {8, 7.5}/n_{5}/east/0/{0,-1.2},       
                {15, 8.25}/n_{6}/north/0/{0,1.2},  
                {12, 4.5}/n_{7}/north west/0/{-1.2,0},
                {15, 6}/n_{8}/north/0/{0,1.2},    
                {8, 10.125}/n_{9}/east/0/{1.2,0},        
                {12, 10.125}/n_{10}/east/0/{1.2,0},     
                {8, 12}/n_{11}/east/0/{1.2,0},
                {12, 12}/n_{12}/east/0/{1.2,0},
                {15, 14.25}/n_{13}/north/0/{0,1.2},
                {6, 14.25}/n_{14}/south/0/{0,-1.2},
                {3, 14.25}/n_{15}/south east/0/{1.2,0},
                {4.5, 16.5}/n_{16}/west/0/{-1.2,0}, 
                {1.5, 18.75}/n_{17}/east/0/{0,-1.2},
                {4, 20.25}/n_{18}/east/0/{1.2,0},
                {8, 16.5}/n_{19}/east/0/{1.2,0}, 
                {12, 16.5}/n_{20}/east/0/{1.2,0}, 
                {8, 20.25}/n_{21}/east/0/{1.2,0},
                {12, 20.25}/n_{22}/east/0/{1.2,0},
                {15, 18}/n_{23}/south/0/{0,-1.2},
                {3, 12}/n_{24}/north east/0/{1.2,0}}
            {
                \node[anchor=\direction] (\name) at (\start) {};
                \draw [ultra thick] (\start) -- ($(\start)+(\shift)$);
            }


    \end{circuitikz}
    
    
} 
                \caption{NCP-KMeans}
                \label{fig:NP_NCP_kmeans}
            \end{subfigure}
            \begin{subfigure}{0.19 \textwidth}
            \centering
                \resizebox{.9\columnwidth}{!}{%

     \begin{circuitikz}  

            \fill [IEE_blue!40!white] (2.5 + 15,4) -- (2.5 + 15, 10.5) -- (9 + 15, 13.124) -- (15.5 + 15, 13.124) -- (15.5 + 15, 7.75) -- (15.5 + 15, 7.75)  -- (15.5 + 15,5.5) -- (11 + 15, 5.5) -- (8.5 + 15,4) -- (2.5 + 15,4);
            
            \fill [IEE_red!40!white]  (17.5,10.7) -- (19.7, 11.6) -- (19.7, 14.5) -- (21, 16.02) -- (30.5, 16.02) -- (30.5, 20.75) -- (17.5, 20.75) -- (17.5, 20) -- (16, 19.15) -- (16, 17.15) -- (17.5, 16) -- (17.5, 13.225) ;
            
            \fill [IEE_green!40!white]  (27, 15.82) -- (30.5, 15.82) -- (30.5, 13.324) -- (27, 13.324) ;
            
            \fill [IEE_orange!40!white] (4.9 + 15, 11.681) -- (4.9 + 15, 14.4) -- (6.15 + 15, 15.82) -- (7 + 15, 15.82) -- (7 + 15, 12.54) -- (4.9 + 15, 11.681)   ;
            
            \fill [IEE_light_blue!40!white] (23.85, 4) -- (26.05, 5.3) -- (30.5, 5.3) -- (30.5, 4) --  (23.85, 4)  ;

            \node[draw, circle, inner sep=2pt] at (14.5 + 15 ,19.5) {\Huge 1};
            \node[draw, circle, inner sep=2pt] at (13   + 15 , 15) {\Huge 2};
            \node[draw, circle, inner sep=2pt] at (6.3 + 15 ,15) {\Huge 3};
            \node[draw, circle, inner sep=2pt] at (14.5 + 15 ,12) {\Huge 4};
            \node[draw, circle, inner sep=2pt] at (14.5 + 15 ,4.65) {\Huge 5};

            \draw (3   + 15 + 0.9, 4.5      )  -- +(0, 0.25)   -- (8    + 15 - 0.9 , 4.5+0.25)        -- +(0, -0.25);    
            \draw  (3   + 15 + 0.3, 10.125   )  -- +(0, -0.25)  -- (3    + 15 + 0.3 , 4.5+0.25)        -- +(0, -0.25);     
            \draw (8   + 15 - 0.6, 4.5      )  -- +(0, 0.25)   -- (4.25 + 15       , 7.5-0.9)          -- +(-0.25, 0);     
            \draw (3   + 15 + 0.6, 4.5      )  -- +(0, 0.25)   -- (8    + 15 - 0.25, 7.5-0.9)        -- +(0.25, 0);      
            \draw (8   + 15 - 0.3, 4.5      )  -- +(0, 0.25)   -- (15   + 15 - 0.25, 8.25+0.3)      -- +(0.25, 0);   
            \draw (3   + 15 + 0.9, 10.125   )  -- +(0, -0.25)  -- (8    + 15 + 0.3 , 10.125-0.25)     -- +(0, 0.25);     
            \draw (4   + 15 - 0.4, 10.125   )                  -- (3    + 15 + 0.6 , 12);                                 
            \draw (4   + 15      , 7.5-0.3  )  -- +(0.25, 0)   -- (8    + 15 + 0.6 , 10.125-0.25)     -- +(0, 0.25);       
            \draw (8   + 15      , 7.5-0.3  )  -- +(0.25, 0)   -- (12   + 15 + 0.3 , 10.125-0.25)    -- +(0, 0.25);    
            \draw [color=red] (15  + 15      , 8.25+0.9 )  -- +(-0.25, 0)  -- (12   + 15 + 0.9 , 10.125-0.25)    -- +(0, 0.25); 
            \draw [color=red](12  + 15 - 0.3, 4.5      )  -- +(0, 0.25)   -- (15   + 15 - 0.25, 6+0.3)         -- +(0.25, 0); 
            \draw (15  + 15      , 6+0.6    )  -- +(-0.25, 0)  -- (8    + 15 + 0.9 , 10.125-0.25)     -- +(0, 0.25);   
            \draw (15  + 15      , 6+0.9    )  -- +(-0.25, 0)  -- (12   + 15 + 0.6 , 10.125-0.25)    -- +(0, 0.25); 
            \draw (8   + 15 + 0.3, 10.125   )                  -- (8    + 15 +0.3, 12);  
            \draw (8   + 15 + 0.9, 10.125   )  -- +(0, 0.25)   -- (12   + 15 + 0.3 , 12-0.25)        -- +(0, 0.25);      
            \draw (12  + 15 + 0.3, 10.125   )  -- +(0, 0.25)   -- (8    + 15 + 0.9 , 12-0.25)         -- +(0, 0.25);      
            \draw (12  + 15 + 0.9, 10.125   )                  -- (12   + 15 + 0.9 , 12);                                
            \draw (8   + 15 + 0.9, 12       )  -- +(0, 0.25)   -- (15   + 15 - 0.25, 14.25+0.6)     -- +(0.25, 0);      
            \draw (6   + 15      , 14.25-0.9)  -- +(0.25, 0)   -- (8    + 15 + 0.3 , 12+0.25)         -- +(0, -0.25);     
            \draw (12  + 15 + 0.9, 12       )  -- +(0, 0.25)   -- (15   + 15 - 0.25, 14.25+0.3)     -- +(0.25, 0);      
            \draw (12  + 15 + 0.6, 12       )  -- +(0, 0.25)   -- (15   + 15 - 0.25, 18-0.6)        -- +(0.25, 0);      
            \draw [color=red](15  + 15      , 14.25+0.9)  -- +(-0.25, 0)  -- (15   + 15 - 0.25, 18-0.9)        -- +(0.25, 0);      
            \draw [color=red](6   + 15      , 14.25-0.3)  -- +(-0.25, 0)  -- (4.5  + 15 - 0.3 , 16.5-0.25)     -- +(0, 0.25);      
            \draw (3   + 15 + 0.6, 14.25    )  -- +(0, 0.25)   -- (4.5  + 15 - 0.9 , 16.5-0.25)     -- +(0, 0.25); 
            \draw [color=red] (3   + 15 + 0.9, 14.25    )  -- +(0, 0.25)   -- (8    + 15 + 0.6 , 20.25-0.25)      -- +(0, 0.25); 
            \draw (3   + 15 + 0.6, 14.25    )                  -- (3    + 15 + 0.6 , 12);                            
            \draw (1.5 + 15      , 18.75-0.9)  -- +(0.25, 0)   -- (4.5  + 15 - 0.9 , 16.5+0.25)     -- +(0, -0.25);  
            \draw (4.5 + 15 - 0.3, 16.5     )  -- +(0, 0.25)   -- (8    + 15 + 0.3 , 16.5+0.25)       -- +(0, -0.25);  
            \draw (1.5 + 15      , 18.75-0.3)  -- +(0.25, 0)   -- (4    + 15 + 0.3 , 20.25-0.25)      -- +(0, 0.25);  
            \draw (1.5 + 15      , 18.75-0.6)  -- +(0.25, 0)   -- (12   + 15 + 0.9 , 20.25-0.25)     -- +(0, 0.25);  
            \draw (4   + 15 + 0.9, 20.25    )  -- +(0, -0.25)  -- (8    + 15 + 0.3 , 20.25-0.25)      -- +(0, 0.25);  
            \draw (8   + 15 + 0.9, 16.5     )  -- +(0, 0.25)   -- (12   + 15 + 0.3 , 16.5+0.25)      -- +(0, -0.25); 
            \draw (12  + 15 + 0.9, 16.5     )  -- +(0, 0.25)   -- (15   + 15 - 0.25, 18-0.3)        -- +(0.25, 0);  
            \draw (8   + 15 + 0.9, 20.25    )  -- +(0, -0.25)  -- (12   + 15 + 0.3 , 20.25-0.25)     -- +(0, 0.25); 

                
            \foreach \start/\name/\direction/\rot/\shift  in {
                {  3 + 15 , 4.5}/n_{1}/north east/0/{1.2,0},
                {  8 + 15 , 4.5}/n_{2}/north west/0/{-1.2,0},
                {  3 + 15 , 10.125}/n_{3}/east/0/{1.2,0},       
                {  4 + 15 , 7.5}/n_{4}/east/0/{0,-1.2},      
                {  8 + 15 , 7.5}/n_{5}/east/0/{0,-1.2},       
                { 15 + 15 , 8.25}/n_{6}/north/0/{0,1.2},  
                { 12 + 15 , 4.5}/n_{7}/north west/0/{-1.2,0},
                { 15 + 15 , 6}/n_{8}/north/0/{0,1.2},    
                {  8 + 15 , 10.125}/n_{9}/east/0/{1.2,0},        
                { 12 + 15 , 10.125}/n_{10}/east/0/{1.2,0},     
                {8   + 15 , 12}/n_{11}/east/0/{1.2,0},
                {12  + 15 , 12}/n_{12}/east/0/{1.2,0},
                {15  + 15 , 14.25}/n_{13}/north/0/{0,1.2},
                {6   + 15 , 14.25}/n_{14}/south/0/{0,-1.2},
                {3   + 15 , 14.25}/n_{15}/south east/0/{1.2,0},
                {4.5 + 15 , 16.5}/n_{16}/west/0/{-1.2,0}, 
                {1.5 + 15 , 18.75}/n_{17}/east/0/{0,-1.2},
                {4   + 15 , 20.25}/n_{18}/east/0/{1.2,0},
                {8   + 15 , 16.5}/n_{19}/east/0/{1.2,0}, 
                {12  + 15 , 16.5}/n_{20}/east/0/{1.2,0}, 
                {8   + 15 , 20.25}/n_{21}/east/0/{1.2,0},
                {12  + 15 , 20.25}/n_{22}/east/0/{1.2,0},
                {15  + 15 , 18}/n_{23}/south/0/{0,-1.2},
                {3   + 15 , 12}/n_{24}/north east/0/{1.2,0}}
            {
                \node[anchor=\direction] (\name) at (\start) {};
                \draw [ultra thick] (\start) -- ($(\start)+(\shift)$);
            }
        

    \end{circuitikz}
    
    
} 
                \caption{NCP-ANAC}
                \label{fig:NP_NCP_AHC}
            \end{subfigure}
            \caption{Grid partitionings of the IEEE RTS 24-Bus System resulting from the aggregation methods stated in table~\ref{tab:clustering_methods} for $\tilde{N} = 5$ nodes in the aggregated grid. Congested lines are marked in red.}
            \label{fig:network_partitioning_comparison}
        \end{figure*}
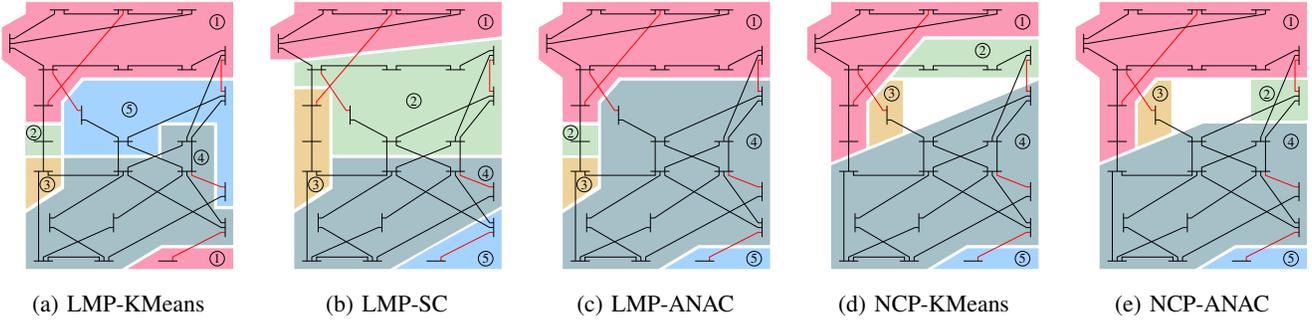
        Based on the introduced IEEE RTS 24-Bus System, we now compare the aggregation methods presented in Section~\ref{s:clustering_methods}.
        First, to illustrate the technical disparities between the aggregation methods, we show the resulting grid partitions for~$\tilde{N} = 5$ nodes in Figure~\ref{fig:network_partitioning_comparison}. We choose this number of nodes to illustrate the different characteristics of the grid partitioning methods while maintaining a clear visual representation. 
        We first look into the partitionings obtained using the LMP metric in Figures~\ref{fig:NP_lmps_kmeans}--\ref{fig:NP_lmps_AHC}. For LMP-KMeans, which does not consider the underlying network topology, we observe that node~$n_7$ (in the bottom right) is grouped with nodes at the top, despite having no direct connection (similar to the nodes in Group~5). This is because their LMPs are similar, as illustrated in Figure~\ref{fig:large_case}. On the other hand, the aggregation method based on Spectral clustering, LMP-SC, shown in Figure~\ref{fig:NP_lmps_spectral}, does not group nodes which are not directly connected, but seems to be less congestion sensitive and solely retains one congested line in the aggregated model. However, LMP-ANAC also enforces a direct connection between nodes within a group while retaining more line congestions than the Spectral clustering. From an electrical perspective, allocating nodes without a direct connection into a single group corresponds to a short circuit between them. This prohibits clear physical interpretations and may challenge the reasonability of group assignments.
        
        Looking at the partitionings obtained with the novel NCP metric in Figures~\ref{fig:NP_NCP_kmeans} and \ref{fig:NP_NCP_AHC}, we observe several similarities for this number of groups. Note, again, that, despite not being shown here, the NCP-KMeans method does not necessarily enforce a direct connection of nodes within a group, although the node clusters are located in similar regions of the grid. This can be observed for other reduction degrees~$\tilde{N}$. Because compared to the Spectral clustering, where the topology information is solely based on the incidence matrix, uses the NCP the electrical distance from the PTDF matrix as an affinity attribute. 
        Comparing the assignments between the LMP- and NCP-based methods, we observe that the latter conserves congested lines better, such that these still exist in the aggregated grid. For example, the LMP-ANAC method, whose results appear most similar to those using the NCP metric, allocates nodes~$n_3$ and~$n_{24}$ into individual groups despite connecting lines having comparable low line loading, as shown in Figure~\ref{fig:large_case}. The NCP metric, therefore, appears to be more congestion-sensitive, which leads to a more similar dispatch between the FM and the AM~\cite{cao_incorporating_2018}.

    \subsection{Results: Aggregated Model Performance}
    \label{s:results_grid_aggregation}
\begin{figure}[t]
            \vspace{-0.8cm} 
            \centerline{
            \includegraphics[width=0.45\textwidth]{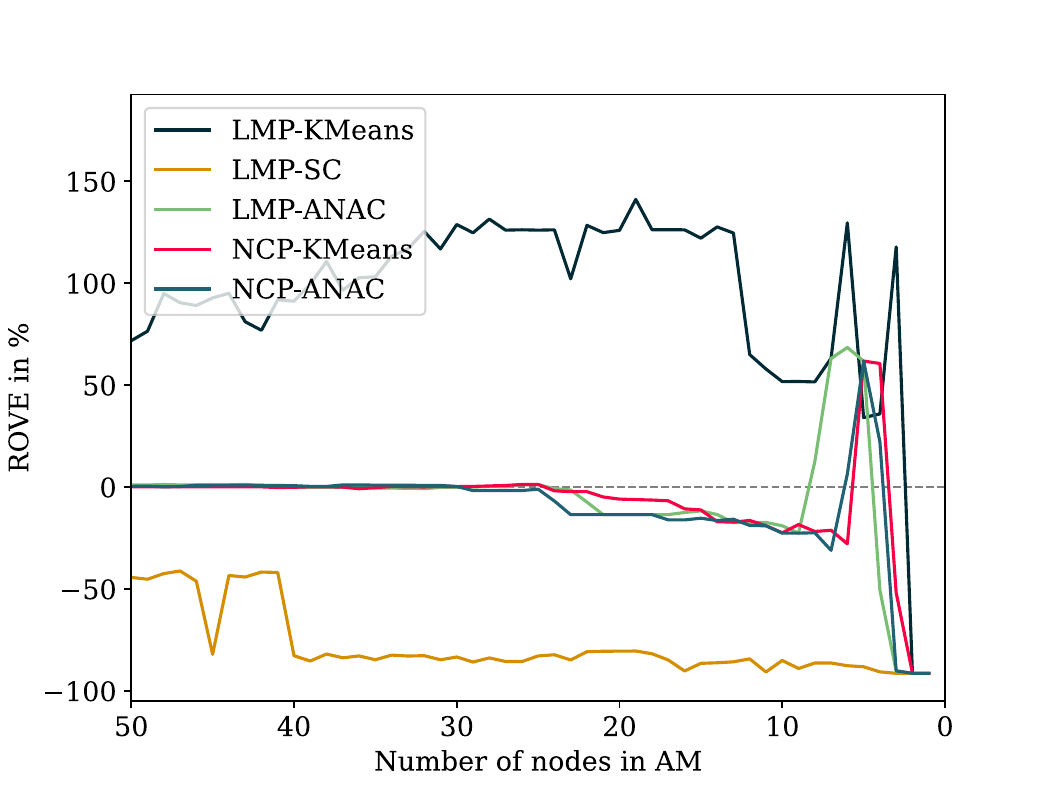}}
            \caption{Relative OFV error for different grid aggregation methods applied to the modified IEEE 300-Bus System as a function of the number of nodes~$\tilde N=50,\dots,1$ in the AM.}
            \label{fig:OFV_rel_error}
            \vspace{-0.55cm}
        \end{figure}

        We assess the quality of the aggregation methods in Table~\ref{tab:clustering_methods} by comparing the FM and AM results for the modified IEEE 300-Bus system with 17 line congestions based on the error metrics ROVE and MRLLV defined in Equations~\eqref{eq:ovf_relative_error} and~\eqref{eq:powerFlow_relative_error} in Section~\ref{s:evaluation}, respectively. Figure~\ref{fig:OFV_rel_error} illustrates the ROVE as a function of the number of nodes in the AMs, for $\tilde{N} < 50$. The strong deviations of the ROVE for high levels of reduction~($\tilde N<7$) indicate a minimum number of nodes needed to represent the grid adequately. Below that, a meaningful aggregation is hardly possible. We observe that the LMP-KMeans tends to overestimate, while the LMP-SC underestimates the system operation cost. However, the ANAC algorithm for both distance metrics and the NCP-KMeans tends to underestimate the cost slightly between $25 \geq \tilde{N} \geq 10$.
        The ROVE results further demonstrate that the LMP-KMeans method produces aggregations that inadequately represent the FM having multiple line congestions. Especially for the aggregations from~\num{20} down to~\num{13} nodes, exhibiting a ROVE of~\SI{140}{\percent}.
        The LMP-SC method also leads to a very poor representation of the FM compared to LMP-KMeans because it attains already for aggregations with fewer nodes than~\num{40} a similar ROVE than the single node case. However, a unit dispatch in the AM with lower total system costs would violate line limit constraints and therefore lead to redispatch cost in a more practical scenario.
        The ANAC method increases the performance compared to the KMeans and Spectral clustering of the LMPs. In this case, both the LMPs and the $\NCP$ demonstrate good performance, showing relatively smooth changes between aggregations.
        Based on the ROVE, the NCP-KMeans shows similar results compared to both ANAC methods. However, the NCP-KMeans method has the advantage of a faster grid partitioning time compared to the ANAC methods, shown in Table~\ref{tab:mean_errors}.
        
\begin{figure}
            \vspace{-0.8cm}
            \centerline{\includegraphics[width=0.45\textwidth]{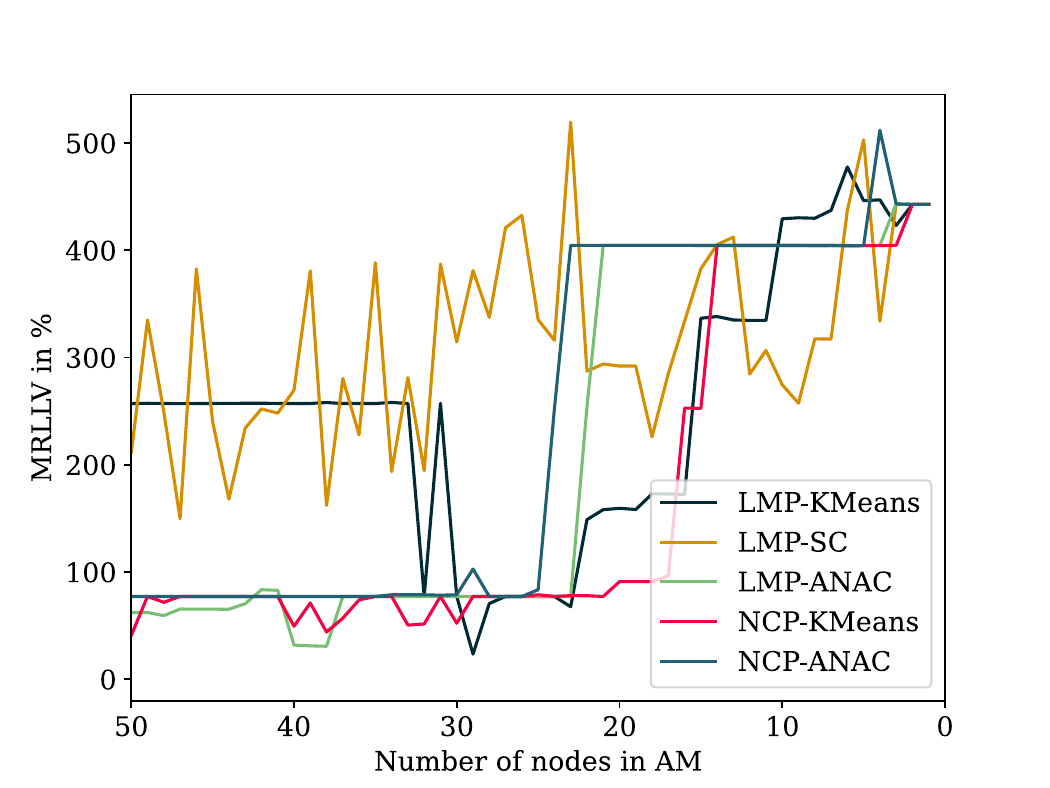}}
            \caption{The maximum relative line limit violation for grid aggregations of the adapted IEEE 300-Bus system, depending on the number of nodes~$\tilde N=50,\dots,1$ in the AM.}
            \label{fig:PowerFlow_rel_error}
            \vspace{-0.5cm}
        \end{figure}

        As illustrated in Figure~\ref{fig:PowerFlow_rel_error}, the MRLLV varies strongly with the different partitioning methods and the number of nodes. The LMP-KMeans exhibits a high violation for a high number of nodes in the AM, then dropping to a similar level than the ANAC methods, before approaching the value of the single node case. The LMP-SC exhibits very high and strongly fluctuating MRLLVs. Solely, the ANAC methods and NCP-KMeans feature MRLLVs with less deviation and lower rates for aggregations with a small level of reduction. Until having a big step-up around~\num{20} nodes in the AM and almost reaching the same MRLLV as the single node case. Which is the same level of reduction, where the ROVE begins to stronger deviate from zero, shown in Figure~\ref{fig:OFV_rel_error}. Notably, already small reductions in the number of nodes result in MRLLV around~\SI{90}{\percent}. This may likely be attributed to the network state of the FM, which already includes several congestions. As a result, the operational schedules of the AMs may lead to physically infeasible power flows on the full grid and may require extensive redispatch measures.
        Table~\ref{tab:mean_errors} summarizes the performance comparison of the presented grid aggregation methods, including the average absolute ROVE and average MRLLV, as well as the average grid partitioning time for the aggregations with~$\tilde N = 50,\dots,1$. 
        As shown previously, the KMeans and Spectral clustering of LMPs also results in the highest average error for both error metrics, whereas the three other clustering methods achieve better performance with less error. The grid partitioning, however, is fast for both KMeans methods as well as the SC method, but very slow for both ANAC methods. The iterative structure of the ANAC algorithm causes a high GPT because it determines every aggregation from $N$ to $\tilde N$, which is not the case for the other partitioning methods. The NCP-KMeans combination, therefore, exhibits a good balance of aggregation performance and partitioning time.
        \begin{table}[!t] \renewcommand{\arraystretch}{1.3} \caption{Absolute relative objective function value error ROVE, maximum relative line limit violation and the grid partitioning time averaged over the aggregations $\tilde{N}=50,\dots,1$ of the IEEE 300-Bus System.} \label{tab:mean_errors} 
        \centering 
        \begin{tabular}{l|c|c|c|c|c} \hline
                \multirow{2}{*}{\textbf{$\,$}} & \multicolumn{3}{c|}{\textbf{LMP}} & \multicolumn{2}{c}{\textbf{NCP}} \\\cline{2-6}
                 & KMeans & SC & ANAC & KMeans & ANAC   \\ \hline
                $\mathrm{ROVE}(z)$ in \%  & \num{93.5}& \num{77.6}& \num{6.4} & \num{6.9}& \num{9.6} \\
                \hline
                 $\mathrm{MRLLV}(\boldsymbol{f})$ in \%& \num{255}& \num{317}& \num{216}& \num{174} & \num{236}  \\ \hline
                 GPT in s & \num{0.33} & \num{0.42} & \num{11.6} & \num{0.34} & \num{11.60}\\
                 \hline 
            \end{tabular} 
            \vspace{-0.3cm}
        \end{table}
        %

\section{Conclusion}\label{s:conclusion}

    In this paper, we show that using LMPs as a distance metric for grid aggregation for DC-OPF may fail to adequately preserve essential physical characteristics of the network topology when there is more than one line congestion. 
    To overcome this problem, we propose the novel Network Congestion Price matrix as a metric for grid partitioning for DC-OPF. Based on two illustrative case studies, we demonstrate the superior performance of the $\NCP$ metric compared to LMPs in terms of relative error in the objective function value and maximum line limit violations. The $\NCP$ metric also is found to be more congestion sensitive than Spectral clustering and achieves a faster partitioning time than the ANAC method.
    Future research should extend the $\NCP$ methodology for multiple time steps and further investigate the robustness when considering multiple periods with varying operational points. Additionally, it should investigate an injection weighted reduced PTDF matrix for better consideration of diverging nodal injections in combination with multiple time steps. Furthermore, guided by the idea of preserving the impact of individual line congestions, a distance metric solely based on a priori available information should be considered for cases where solutions for the DC-OPF problem are unavailable, e.g., in planning problems.

\appendices

    \section{Grid Aggregation \& Equivalent Susceptances}\label{s:appendix_grid_aggregation}
    
        Network reduction methods, including Ward- and REI-Reduction, have been utilized for considerable time~\cite{ward_equivalent_1949, debs_estimation_1975}. These methods simplify the power flow physics adapting solely the parameters of the lines connecting the nodes to be aggregated, but not those of the surrounding grid.
        This may be adequate when reducing a grid outside the area of interest. However, when considering the complexity of the full grid, this can result in inappropriate grid representations~\cite{fortenbacher_transmission_2018}.
        To overcome this issue, grid aggregations for DC-OPF should be done based on reducing the PTDF matrix, as shown in ~\eqref{eq:reduced_PTDF}, which considers the impact of a nodal power injection on all lines. 
        Combined with determining equivalent susceptances for bus-angle DC-OPF models for faster computation~\cite{shi_novel_2015, fortenbacher_transmission_2018}.

    \section{Spectral Clustering}
    \label{sec:appendix_SC}

        For the Spectral Clustering, we implement the algorithm that the authors in~\cite{cao_incorporating_2018} proposed. The algorithm requires a Laplacian matrix~$\mathscr{L} \in \mathbb{R}^{N \times N}$ representing the affinity of the network nodes, where higher values represent a higher affinity and vice versa. We define the inverse LMP difference element wise per line~$\boldsymbol{\rho} \in \mathbb{R}^{L \times 1}$ by:
        \begin{equation}
            {\rho}_l = \big |{\big ( \mathbf{K}\, \mathbf{LMP}\big )_{l}\big |}^{-1}, \quad \forall l \in \mathcal{L} .
        \end{equation}
        We then define the Laplacian matrix~$\mathscr{L}$ using the diagonal matrix of~$\boldsymbol{\rho}$ and the line-node incidence matrix~$\mathbf{K}$ as:
        \begin{equation}\label{eq:sc_laplacian}
            \mathscr{L} = \mathbf{K}\tran \, \operatorname{diag}(\boldsymbol{\rho})\, \mathbf{K}.
        \end{equation}
        For the clustering process, we pass the Laplacian matrix to the spectral clustering algorithm implemented in the Scikit-learn~\cite{pedregosa_scikit-learn_nodate} with the \textit{affinity} option set to \textit{'precomputed'}. The algorithm normalizes the Laplacian matrix before applying a KMeans algorithm to the $\tilde N$-smallest eigenvalues of~$\mathscr{L}$.

    \section{Adjacent Node Agglomerative Clustering}
    \label{sec:Appendix_ACA}
        In the following, we provide a detailed explanation of the Adjacent Node Agglomerative Clustering given in Algorithm~\ref{alg:anac} below. Let~$c$ be the cluster counter, starting with the number of nodes in the full grid~$N$ and decreasing every iteration until it reaches the desired number of node clusters in the aggregated grid~$\tilde N$. Secondly, we define~$\mathbf{i}_c$ as a vector of ones with size~$N$, denoted as~$\mathbf{e}^N \in \mathbb{Z}^{N\times 1}$. Additionally, we define~$\mathbf{X}_c$ either by the $\mathbf{LMP} \in \mathbb{R}^{N \times 1}$ or the~$\NCP \in \mathbb{R}^{N \times L}$ to calculate the distance~$\boldsymbol{\mathcal{D}} \in \mathbb{R}^{N\times N}$ between each pair of nodes for all entries, denoted by the colon as index.  To ensure that nodes within each cluster have a direct line connection, we define a weighted distance matrix~$\mathbf{W} \in \mathbb{Z}^{N \times N}$ using the node-incidence matrix~$\mathbf{K}_c^\mathrm{N} = \mathbf{K}_c\tran \, \mathbf{K}_c$, which has entry 1 if two nodes are connected by a line and zero if they have no connection. Intuitively, in each iteration, the Adjacent Node Agglomerative Clustering algorithm groups the two nodes $(\tilde n, \tilde m)$ with the minimum weighted distance~$\mathbf{W}$. Afterward, the LMPs or NCPs~$\mathbf{X}_c$ of the new group are updated using the weighted cluster average. Note, that we do not rerun the DC-OPF between iterations to obtain those values but always use the solutions of the FM.
        \begin{algorithm}
            \caption{Adjacent Node Agglomerative Clustering}\label{alg:anac}
            \begin{algorithmic}[1]
                \State $c \gets N$ \Comment{Cluster counter}
                \State $ \mathbf{i}_c \gets \boldsymbol{e}^{N} $ \Comment{Cluster-size counter}
                \State $ \mathbf{X}_c \gets\mathbf{LMP}$ or $\mathbf{NCP}$ \Comment{Distance metric}
                \State $\mathbf{K}_c^\mathrm{N} \gets \mathbf{K}_c\tran \, \mathbf{K}_c$ 
                \Comment{Initial nodal adjacency}
                \State $\mathbf{M}^\mathrm{nc}_k \gets \mathbf{I}^N$ \Comment{Initialize node mapping}
                \While{$c > \tilde N$}
                    \Statex Determine weighted distance matrix~$\mathbf{W}_c$
                    \State ${\mathcal{D}}_{n,m} \gets \lVert \mathbf{X}_{n,:,c} - \mathbf{X}_{m,:,c} \rVert_2, \quad \forall n,m = 1,\dots, c$
                    \State $\mathbf{W} \gets \boldsymbol{\mathcal{D}} + (\boldsymbol{e}^{c} \, (\boldsymbol{e}^{c})\tran - \mathbf{K}_c^\mathrm{N} + \mathbf{I}^{c}) \, \varepsilon, \quad \varepsilon \geq \operatorname{max}(\boldsymbol{\mathcal{D}})$
                    \Statex Get indices~$(\tilde n, \tilde m)$ of minimum entries of~$\mathbf{W}$
                    \State $(\tilde n,\tilde m) \gets \operatorname{arg}\, \min_{(n,m)} \mathbf{W} \in \mathbb{R}^{c \times c}$
                    \Statex Define node-cluster mapping matrix
                    \State $\mathbf{M}^\mathrm{nc} \gets \bigl[\mathbf{I}^{c-1}_{:,1:\tilde m-1}\, | \, \boldsymbol{0} \, | \, \mathbf{I}^{c-1}_{:,\tilde m:c-1}\bigr]$; $\mathrm{M}^\mathrm{nc}_{\tilde n, \tilde m} \gets 1$
                    \Statex Compute nodal incidence matrix
                    \State $\mathbf{K}^\mathrm{N}_{c-1} \gets \mathbf{M}^\mathrm{nc} \, \mathbf{K}_{c}^\mathrm{N} \, (\mathbf{M}^\mathrm{nc})\tran$
                    \Statex Determine cluster size counter and cluster weights
                    \State $\mathbf{i}_{c-1} \gets \mathbf{M}^\mathrm{nc} \, \mathbf{i}_{c}$
                    \State $\mathbf{w} \gets \boldsymbol{e}^{c}$; $\mathrm{w}_{\tilde{m}} \gets 1/\mathrm{i}_{\tilde{n},c-1}$; $\mathrm{w}_{\tilde{n}} \gets 1-\mathrm{w}_{\tilde{m}}$
                    \Statex Compute weighted average of the distance metric
                    \State $\mathbf{X}_{c-1} \gets \mathbf{M}^\mathrm{nc}\, \mathbf{X}_{c} \, \operatorname{diag}(\mathbf{w})$
                    \State $\mathbf{M}^\mathrm{nc}_k \gets \mathbf{M}^\mathrm{nc}_k \, \mathbf{M}^\mathrm{nc}$
                    \State $c \gets c-1$
                \EndWhile
                \State \textbf{return} $\mathbf{M}^\mathrm{nc}_k$ 
            \end{algorithmic}
        \end{algorithm}

    \section*{References}
    \addcontentsline{toc}{section}{References}
    
    \printbibliography[heading=none]
    
\end{document}